\pdfoutput=1 
\documentclass{amsart}%
\usepackage{amsfonts}
\usepackage{graphicx}
\usepackage{amscd}
\usepackage{amsmath}
\usepackage{amssymb}%
\providecommand{\U}[1]{\protect\rule{.1in}{.1in}}

\theoremstyle{plain}

\newtheorem{remark}{Remark}

\numberwithin{equation}{section}

\ifx\pdfoutput\relax\let\pdfoutput=\undefined\fi
\newcount\msipdfoutput
\ifx\pdfoutput\undefined\else
\ifcase\pdfoutput\else
\ifx\paperwidth\undefined\else
\ifdim\paperheight=0pt\relax\else\pdfpageheight\paperheight\fi
\ifdim\paperwidth=0pt\relax\else\pdfpagewidth\paperwidth\fi
\fi\fi\fi
\begin{document}
\title[\textsf{THE SUBALGEBRA STRUCTURE OF THE }\textsl{trigintaduonions}]{\textsf{THE BASIC SUBALGEBRA STRUCTURE OF THE CAYLEY-DICKSON ALGEBRA OF
DIMENSION 32 (}\textsl{trigintaduonions}\textsf{)}}
\author{Raoul E. Cawagas, \textsl{et al}*}
\address{Raoul E. Cawagas\\
SciTech R\&D Center\\
Polytechnic University of the Philippines\\
Sta. Mesa, Manila}
\email{raoulec1@\thinspace yahoo.com}
\thanks{\textbf{*Co-Authors:} \textbf{Alexander S. Carrascal, Lincoln A. Bautista,
John P. Sta. Maria, Jackie D. Urrutia, Bernadeth Nobles}}

\begin{abstract}
The Cayley-Dickson algebras $\mathbb{R}$ (real numbers), $\mathbb{C}$ (complex
numbers), $\mathbb{H}$ (quaternions), $\mathbb{O}$ (octonions), $\mathbb{S}$
(sedenions) and $\mathbb{T}$ (trigintaduonions) have attracted the attention
of several mathematicians and theoretical physicists because of their
important applications in both pure mathematics and theoretical physics. This
paper deals with the determination of the \emph{basic subalgebra structure} of
the algebra $\mathbb{T}$ by analyzing the loop $T_{{\small L}}$ of order 64
generated by its 32 basis elements. The analysis shows that $T_{L}$ is a
\emph{non-associative finite invertible loop} (NAFIL) with 373 non-trivial
subloops of orders 32, 16, 8, 4, and 2 all of which are normal. These subloops
generate subalgebras of $\mathbb{T}$ of dimensions 16, 8, 4, 2, and 1 which
form the elements of its structure.

\end{abstract}
\maketitle

\section{Introduction}

The Cayley-Dickson algebras $\mathbb{C}$ (complex numbers 2-D), $\mathbb{H}$
(quaternions 4-D), $\mathbb{O}$ (octonions 8-D), $\mathbb{S}$ (sedenions
16-D), and $\mathbb{T}$ (\emph{trigintaduonions} 32-D)\footnote{The
32-dimensional Cayley-Dickson algebra known as the $trigintaduonions$ (from
the Latin word $trigintaduo,$ meaning 32) is also called the $2^{5}$-$ions.$}
are real algebras obtained from the real numbers $\mathbb{R}$ (1-D) by a
doubling procedure called the Cayley-Dickson (C-D) process [1, 7]. Thus we
have the following C-D doubling chain:

\begin{center}
$\mathbb{R}$\textbf{\ }$\mathbf{\subset}$\textbf{\ }$\mathbb{C}$%
\textbf{\ }$\mathbf{\subset}$\textbf{\ }$\mathbb{H}$\textbf{\ }%
$\mathbf{\subset}$\textbf{\ }$\mathbb{O}$\textbf{\ }$\mathbf{\subset}%
$\textbf{\ }$\mathbb{S}$\textbf{\ }$\mathbf{\subset}$\textbf{\ }$\mathbb{T}%
$\textbf{\ }$\mathbf{\subset}$\textbf{\ }...
\end{center}

\noindent This shows that the trigintaduonions $\mathbb{T}$ contains
$\mathbb{S}$, $\mathbb{O}$, $\mathbb{H}$, $\mathbb{C}$, and $\mathbb{R}$ as
subalgebras. These, however, are not the only subalgebras of $\mathbb{T}$: any
subalgebra of $\mathbb{S}$, $\mathbb{O}$, $\mathbb{H}$, and $\mathbb{C}$ is
also a subalgebra of $\mathbb{T}$ as well as others generated by its basis elements.

In a previous study [2, 3] we determined (using the software \textsl{FINITAS
}[4]) the \emph{basic subalgebra structure} of the C-D sedenion algebra
$\mathbb{S}$ of dimension 16. We showed that it contains an embedded loop
$S_{{\small L}}$ of order 32 (called the \textsl{standard sedenion loop})
generated by its 16 basis elements. Analysis of $S_{L}$ showed that it
contains 15 maximal subloops of order 16 all of which are non-abelian NAFILs
(\textsl{non-associative finite invertible loops}) [5]. Of these, 8 are
isomorphic to the standard \textsl{octonion loop} $O_{L}$ of order 16
generated by the 8 basis elements of $\mathbb{O}$ while 7 are isomorphic to a
newly identified loop $\widetilde{O}_{L}$ which we called the
\textsl{quasi-octonion loop. }Up to isomorphisms,\textsl{\ }these subloops
generate 8-dimensional subalgebras, $\mathbb{O}$ and $\widetilde{\mathbb{O}},$
of the sedenions $\mathbb{S}$.

The next C-D algebra in the chain after $\mathbb{S}$ is the
\textsl{trigintaduonion} algebra $\mathbb{T}$ of dimension $32.$ The
captivating thing about $\mathbb{T}$ is that all of the well-known real
division algebras $\mathbb{R}$, $\mathbb{C}$, $\mathbb{H}$, and $\mathbb{O}$
fit nicely inside it as subalgebras. Hence any object involving these algebras
can be dealt with in $\mathbb{T}.$ Moreover, it also contains the sedenions
$\mathbb{S}$ and several sedenion-type algebras that have potential
applications in pure mathematics and theoretical physics. This is our
motivation for studying these algebras.

\section{Methodology}

The important question related to the study of any algebra is to know its
subalgebra composition and lattice which determine its structure.

This paper presents an attempt to determine the \emph{basic subalgebra
structure} of $\mathbb{T}$ by determining the subloop composition and lattice
of the loop $T_{{\small L}}$ of order $64$ generated by the $32$ basis
elements of $\mathbb{T}$. The subloops of $T_{{\small L}}$ are of orders 32,
16, 8, 4, and 2. These generate subalgebras of $\mathbb{T}$ of dimensions 16,
8, 4, 2, and 1, respectively. Thus there is a one-to-one
correspondence\footnote{This correspondnece can be considered as a form of
\emph{duality}.} between the subloops of $T_{{\small L}}$ and the basic
subalgebras of $\mathbb{T}.$ Because of this the basic subalgebra composition
of $\mathbb{T}$ is directly related to that of $T_{{\small L}}.$ To complete
the description of the structure of $T_{{\small L}}$ we must show how its
subloops fit inside it and how they are related to each other by determining
its subloop lattice; and hence, by the correspondence, also the basic
subalgebra lattice of $\mathbb{T}.$

\subsection{Definitions}

In the standard literature the term \emph{algebra }[1] is taken to mean a
finite dimensional vector space over a field \emph{F} (like $\mathbb{R}$) with
a bilinear multiplication (not necessarily associative) and with a unit
element. Being a vector space, we can choose a basis in terms of which each
element of the algebra can be written as a linear combination of the basis
elements. Such an algebra is completely defined by the multiplication rules of
its basis by means of a \emph{multiplication table}. A \emph{loop }[5], on the
other hand, is a binary system satisfying all group axioms but not necessarily
the associative property. For finite loops of small order, multiplication is
usually defined by a \emph{Cayley table}.

The Cayley-Dickson algebras [1] is a sequence of algebras over the real
numbers $\mathbb{R},$ each with twice the dimension of the previous one as
indicated by the doubling chain in Section 1. Thus, these algebras have
dimensions of the form $2^{n},$ where $n\geq0$ is an integer. Any algebra in
the chain contains all of the algebras (and their subalgebras) before it as
subalgebras. For $n\geq3$ such an algebra is \emph{not associative}.

\subsection{Notation}

In what follows, let the $32=2^{5}$ basis elements of the
\emph{trigintaduonion algebra} $\mathbb{T}$ be represented by the set
$T_{{\small E}}=\{e_{0},e_{1},e_{2},e_{3},...,e_{31}\}$, where $e_{0}$ is the
\emph{unit} element. If $e_{i},e_{j},e_{k}\in T_{{\small E}},$ then $e_{i}\ast
e_{j}=\pm e_{k},$ where $-e_{k}\notin T_{{\small E}}$. Because of this the set
$T_{{\small E}}$ is not closed under the operation $\ast$ of trigintaduonion
multiplication since it contains only positive elements. Thus the 32 elements
of $T_{{\small E}}$ must be extended to include 32 negative elements
$\in\mathbb{T}$ to form the \emph{trigintaduonion loop} $T_{{\small L}}$ of
order $64=2^{6}$. Accordingly, these 64 loop elements will be represented by
the set $T_{{\small L}}=\pm\{e_{0},e_{1},e_{2},e_{3},...,e_{31}\}.$ If $i=0,$
then $e_{0}=\mathfrak{1}$ (the \emph{unit}) while if $i\geq1,$ then $\pm
e_{i}$ is an \emph{imaginary} such that $(\pm e_{i})^{2}=-1.$ Moreover,
$(-e_{i})\ast(e_{j})=e_{i}\ast(-e_{j})=-(e_{i}\ast e_{j}).$

\subsection{Computer Use}

In this study we relied heavily on the use of the computer software
\emph{FINITAS} [4], the LOOPS package for GAP [9], and other programs for the
construction and analysis of the various algebraic structures involved.

First, the Cayley table of the trigintaduonion loop $T_{{\small L}}$ was
computer generated by means of a special computer program based on the
Cayley-Dickson process. The result is shown in Table 1 which only shows its
main portion that corresponds to the multiplication table of the 32 basis
elements of $\mathbb{T}.$

Next, using the software \textsl{FINITAS} we determined the subloop
composition of the trigintaduonion loop $T_{{\small L}}$ by generating the
Cayley tables of its subloops. We also analyzed $T_{{\small L}}$ to determine
its basic structural properties (Section 3.5). Then we identified its subloops
and classified them into isomorphy classes by subjecting them to isomorphism
and other tests. For this, we used both FINITAS and LOOPS\footnote{Prof.
Michael Kinyon (Denver University, USA) used the LOOPS package to assist us in
determining the isomorphy classes of the 31 sedenion-type subloops of
$T_{{\small L}}.$}. The data obtained from the various tests conducted were
then analyzed.

\section{Results and Discussions}

The trigintaduonion algebra $\mathbb{T}$ contains an embedded loop
$T_{{\small L}}$ of order 64 generated by its 32 basis elements. Analysis
(using \emph{FINITAS}) shows that this loop is a NAFIL (\emph{non-associative
finite invertible loop}) with 373 non-trivial subloops of orders $m=32$ (31
loops), $m=16$ (155 loops), $m=8$ (155 loops), $m=4$ (31 loops), and $m=2$ (1
loop). The subloops of orders 32 and 16 are NAFILs while those of orders 8, 4,
and 2 are groups. Moreover, all of these 373 subloops of $T_{{\small L}}$ are
\emph{normal}. Such a loop, called a \emph{Hamiltonian} loop, is known to have
a \emph{modular} lattice [6].

\begin{center}%
{\includegraphics[
natheight=2.923900in,
natwidth=4.545500in,
height=2.9239in,
width=4.5455in
]%
{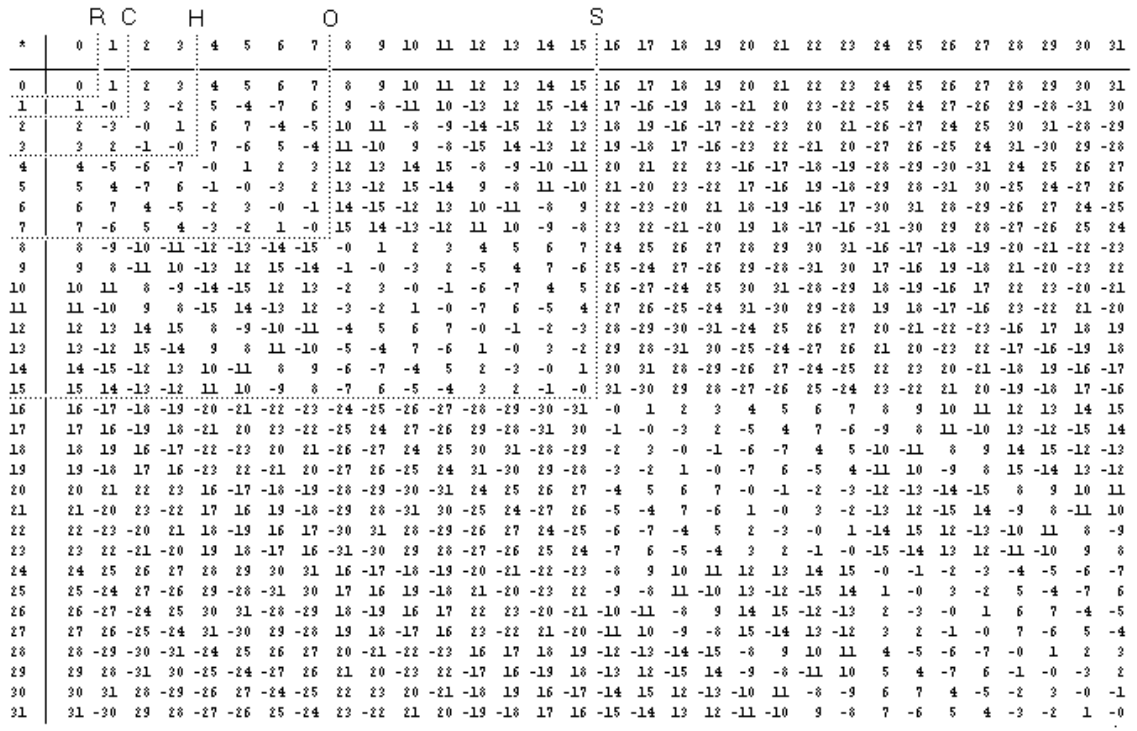}%
}%

\end{center}

\begin{quotation}
Table 1. Main portion of the Cayley table of the trigintaduonion loop
$T_{L}=\pm\{e_{0},e_{1},e_{2},e_{3},...,e_{31}\}$ of order $m=64.$ This
portion corresponds to the multiplication table of the basis $T_{{\small E}%
}=\{e_{0},e_{1},e_{2},e_{3},...,e_{31}\}$ of $\mathbb{T}.$ Note how
$\mathbb{S}$, $\mathbb{O}$, $\mathbb{H}$, and $\mathbb{R}$ are contained in
$\mathbb{T}$. For convenience of notation, we represent each loop element
$e_{i}$ by its subscript $i,$ that is, we set $i=e_{i}.$
\end{quotation}

\subsection{Subloops of Order 32}

The maximal subloops of $T_{{\small L}}${\small \ }are the 31 subloops of
order 32 (called \emph{sedenion-type} loops). One of these is the "standard"
sedenion loop $S_{{\small L}}$ generated by the basis of the sedenion algebra
$\mathbb{S}$. In addition to $S_{{\small L}}$ three more of these 31
sedenion-type loops of order 32 have been identified as distinct
(non-isomorphic). In terms of the basis elements of $\mathbb{T}$, we
have\footnote{\emph{FINITAS} decomposes the loop $T_{{\small L}}$ into its 373
non-trivial subloops numbered \#2 - \#374; subloop \#1 is the trivial subloop
of order 1.}:\medskip

\begin{itemize}
\item $S_{{\small L}}(\#2)=\pm\{0,1,2,3,4,5,6,7,8,9,10,11,12,13,14,15\}$
\end{itemize}

\qquad(std sedenion loop)

\begin{itemize}
\item $S_{{\small L}}^{{\tiny \alpha}}(\#7)=\pm
\{0,1,2,3,8,9,10,11,20,21,22,23,28,29,30,31\}$
\end{itemize}

\qquad($\alpha$-sedenion loop)

\begin{itemize}
\item $S_{{\small L}}^{{\tiny \beta}}(\#10)=\pm
\{0,1,2,3,12,13,14,15,20,21,22,23,24,25,26,27\}$
\end{itemize}

\qquad($\beta$-sedenion loop)

\begin{itemize}
\item $S_{{\small L}}^{{\tiny \delta}}(\#4)=\pm
\{0,1,2,3,4,5,6,7,24,25,26,27,28,29,30,31\}$
\end{itemize}

\qquad($\gamma$-sedenion loop)

\subsubsection{Isomorphy Classes of Sedenion-type Subloops}

Further analysis shows that these four distinct subloops represent exactly
four \emph{isomorphy classes}:

\begin{itemize}
\item $\mathcal{C}_{\cong}\{S_{L}(\#2)\}:$ $\longrightarrow$16 subloops

\item $\mathcal{C}_{\cong}\{S_{L}^{\alpha}(\#7)\}:$ $\longrightarrow$7 subloops

\item $\mathcal{C}_{\cong}\{S_{L}^{\beta}(\#10)\}:$ $\longrightarrow$7 subloops

\item $\mathcal{C}_{\cong}\{S_{L}^{\gamma}(\#4)\}:$ $\longrightarrow$1 subloop
\end{itemize}

\noindent This means that each of the 31 sedenion-type subloops of
$T_{{\small L}}$ belongs to just one, and only one, of these classes. Up to
isomorphism, it can be shown that the subloops $S_{L},$ $S_{L}^{\alpha},$
$S_{L}^{\beta},$ $S_{L}^{\gamma}$ of $T_{{\small L}}$ generate subalgebras
$\mathbb{S},$ $\mathbb{S}^{\alpha},$ $\mathbb{S}^{\beta},$ $\mathbb{S}%
^{\gamma}$ of $\mathbb{T}$.

The Cayley tables of $S_{{\small L}},$ $S_{{\small L}}^{{\tiny \alpha}},$
$S_{{\small L}}^{{\tiny \beta}},$ and $S_{{\small L}}^{{\tiny \gamma}}$ are
shown in Tables 2, 3, 4, and 5. Like Table 1, these Cayley tables only show
their main portions; these correspond to the multiplication tables of the
subalgebras $\mathbb{S},$ $\mathbb{S}^{\alpha},$ $\mathbb{S}^{\beta},$ and
$\mathbb{S}^{\gamma}$ of $\mathbb{T}$.

\begin{center}%
{\includegraphics[
natheight=1.733100in,
natwidth=3.039000in,
height=1.7331in,
width=3.039in
]%
{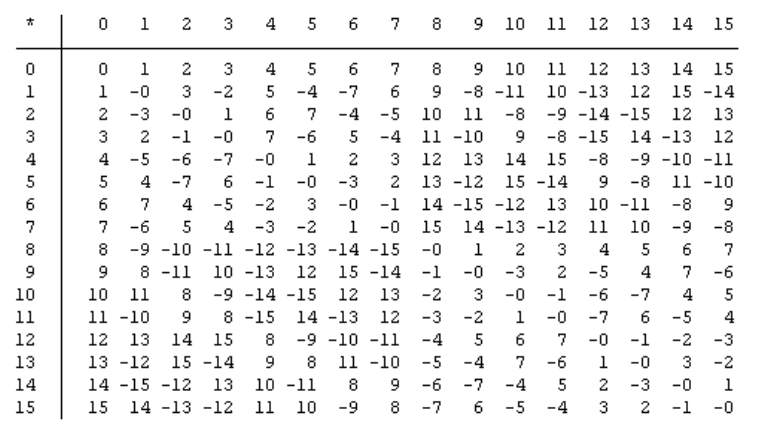}%
}%

\end{center}

\begin{quotation}
Table 2. Main portion of the Cayley table of the standard sedenion loop
$S_{{\small L}}(\#2)=\pm\{0,1,2,3,4,5,6,7,8,9,10,11,12,13,14,15\}$ of order
$n=32.$
\end{quotation}

\begin{center}%
{\includegraphics[
natheight=1.733900in,
natwidth=3.039000in,
height=1.7339in,
width=3.039in
]%
{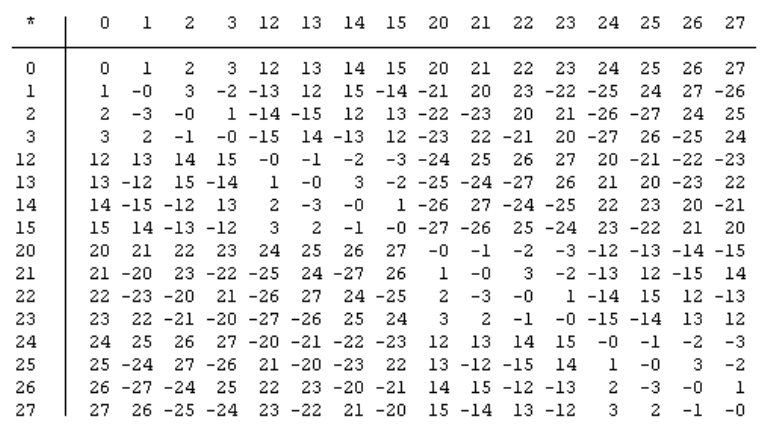}%
}%

\end{center}

\begin{quotation}
Table 3. Main portion of the Cayley table of the $\alpha$-sedenion loop
$S_{L}^{\alpha}(\#7)=\pm\{0,1,2,3,8,9,10,11,20,21,22,23,28,29,30,31\}$ of
order $n=32.$
\end{quotation}

\begin{center}%
{\includegraphics[
natheight=1.733900in,
natwidth=3.039000in,
height=1.7339in,
width=3.039in
]%
{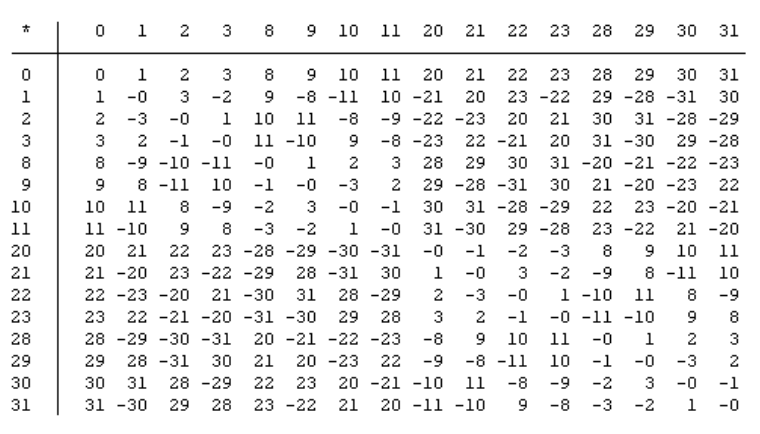}%
}%

\end{center}

\begin{quotation}
Table 4. Main portion of the Cayley table of the $\beta$-sedenion loop
$S_{{\tiny L}}^{{\tiny \beta}}(\#10)=\pm
\{0,1,2,3,12,13,14,15,20,21,22,23,24,25,26,27\}$ of order $n=32.$
\end{quotation}

\begin{center}%
{\includegraphics[
natheight=1.733900in,
natwidth=3.038100in,
height=1.7339in,
width=3.0381in
]%
{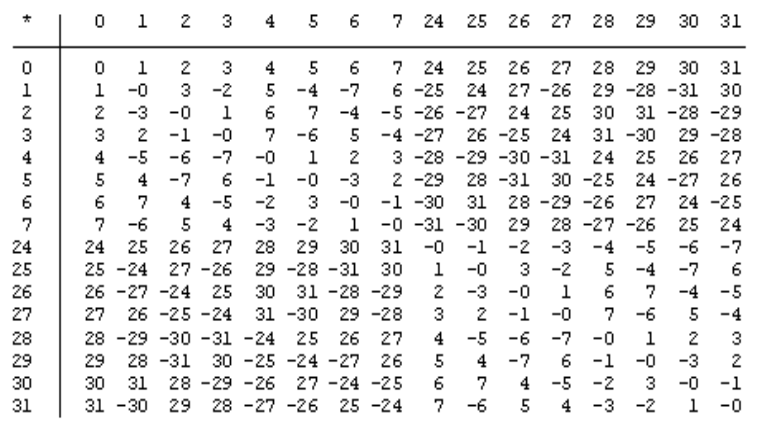}%
}%

\end{center}

\begin{quotation}
Table 5. Main portion of the Cayley table of the $\gamma$-sedenion loop
$S_{{\tiny L}}^{{\tiny \gamma}}(\#4)=\pm
\{0,1,2,3,4,5,6,7,24,25,26,27,28,29,30,31\}$ of order $n=32.$
\end{quotation}

\begin{remark}
The Cayley table of the loop $T_{{\small L}}$ (Table 1) of order $64$ has $64
$ rows and $64$ columns consisting of four portions (partitions) each with 32
rows and 32 columns. Such a table is somewhat large and we only show its main
portion: the multiplication table of the basis $\mathbb{T}_{{\small E}}$ of
$\mathbb{T}$. Similarly, Tables 2, 3, 4, and 5 only show their main portions.
\end{remark}

\subsection{Subloops of order 16}

Analysis shows that $T_{{\small L}}$ has 155 subloops of order 16. These are
\emph{octonion-type} NAFIL loops that form exactly two isomorphy
classes:\bigskip

\begin{itemize}
\item $\mathcal{C}_{\cong}\{O_{{\small L}}\}$ of \emph{octonion} loops [Class
representative: $\#5=\pm\{0,1,2,3,4,5,6,7\}$] = 50 loops

\item $\mathcal{C}_{\cong}\{\widetilde{O}_{{\small L}}\}$ of
\emph{quasi-octonion }loops [Class representative: $\#11=\pm
\{0,1,2,3,12,13,14,15\}$] = 105 loops\bigskip
\end{itemize}

\noindent Thus, every octonion-type loop either belongs to the class
$\mathcal{C}_{\cong}\{O_{{\small L}}\}$ or to the class $\mathcal{C}_{\cong%
}\{\widetilde{O}_{{\small L}}\}.$ These are the maximal subloops of the
sedenion-type loops.\bigskip

\subsection{Subloops of Order 8, 4, and 2}

The subloops of order 8, 4, and 2 are all groups. Analysis shows that all of
the 155 subloops of order $m=8$ are groups isomorphic to the \emph{quaternion
group }$Q_{8}$. On the other hand, the 31 subloops of order $m=4$ are groups
isomorphic to the cyclic group $C_{4};$ the lone subloop of order $m=2$ is a
group isomorphic to the cyclic group $C_{2}$. Thus they form the following
isomorphy classes:

\begin{itemize}
\item $\mathcal{C}_{\cong}\{Q_{{\small 8}}\}$ of quaternion groups
$Q_{{\small 8}}$: $\longrightarrow$155 subloops

\item $\mathcal{C}_{\cong}\{C_{{\small 4}}\}$ of cyclic groups $C_{{\small 4}%
}$: $\longrightarrow$31 subloops

\item $\mathcal{C}_{\cong}\{C_{{\small 2}}\}$ of cyclic group $C_{2}$:
$\longrightarrow$1 subloop
\end{itemize}

\subsection{The Isomorphy Classes of the Subloops of $T_{{\protect\small L}}$}

We have now identified all of the 9 isomorphy classes of the 373 non-trivial
subloops of $T_{{\small L}}:$\medskip

\begin{itemize}
\item $\mathcal{C}_{\cong}\{S_{{\small L}}\},$ $\mathcal{C}_{\cong%
}\{S_{{\small L}}^{\alpha}\},$ $\mathcal{C}_{\cong}\{S_{{\small L}}^{\beta
}\},$ $\mathcal{C}_{\cong}\{S_{{\small L}}^{\gamma}\},$ $\mathcal{C}_{\cong%
}\{O_{{\small L}}\},$ $\mathcal{C}_{\cong}\{\widetilde{O}_{{\small L}}\},$
$\mathcal{C}_{\cong}\{Q_{{\small 8}}\},$ $\mathcal{C}_{\cong}\{C_{{\small 4}%
}\},$ and $\mathcal{C}_{\cong}\{C_{{\small 2}}\}.$
\end{itemize}

\noindent where $S_{{\small L}},$ $S_{{\small L}}^{\alpha},$ $S_{{\small L}%
}^{\beta},$ $S_{{\small L}}^{\gamma},$ $O_{{\small L}},$ $\widetilde
{O}_{{\small L}},$ $Q_{{\small 8}},$ $C_{{\small 4}},$ and $C_{{\small 2}}$
are class representatives. Here, $S_{{\small L}}$ is the standard sedenion
loop of order 32 while $S_{{\small L}}^{\alpha},$ $S_{{\small L}}^{\beta},$
$S_{{\small L}}^{\gamma}$ are newly identified loops of order 32. The loop
$O_{{\small L}}$ of order 16 is the standard octonion loop while
$\widetilde{O}_{{\small L}}$ is called the \emph{quasi-octonion} loop of order
16. On the other hand, $Q_{{\small 8}}$ is the quaternion group of order 8,
$C_{{\small 4}}$ is the cyclic group of order 4, and $C_{{\small 2}}$ is the
cyclic group of order 2.

It can be shown that each of the 373 subloops under these classes generate
subalgebras of $\mathbb{T}$. Thus the basic subalgebra composition of the
algebra $\mathbb{T}$ corresponds to the subloop composition of the loop
$T_{{\small L}}.$ This also means that each of the subloop isomorphy classes
of $T_{{\small L}}$ determines a class of isomorphic subalgebras of
$\mathbb{T}.$

\subsection{Maximal Subloop Compositions of Sedenion-Type Loops}

Analysis also shows that each of the 31 sedenion-type subloops of
$T_{{\small L}}$ belong to exactly one of the following three \emph{maximal
subloop composition types}:\medskip

\begin{itemize}
\item $[$8 $O_{{\small L}}$ + 7 $\widetilde{O}_{{\small L}}]=[8+7]$
:$\longrightarrow$ 17 subloops

\item $[$2 $O_{{\small L}}$ + 13 $\widetilde{O}_{{\small L}}]=[2+13]$
:$\longrightarrow$ 7 subloops

\item $[$0 $O_{{\small L}}$ + 15 $\widetilde{O}_{{\small L}}]=[0+15]$
:$\longrightarrow$ 7 subloops\medskip
\end{itemize}

Thus we find that any sedenion-type loop contains at least 7 quasi-octonion
loops $\widetilde{O}_{{\small L}}$ as subloops. Table 6 shows the four
identified sedenion-type loops in Section 3.1 and their subloop composition
types.\bigskip

\begin{center}%
\begin{tabular}
[c]{|l||c|c|c|c|}\hline
\quad Loop\quad\quad & \quad$S_{{\small L}}(\#2)$\quad\quad & \quad
$S_{{\small L}}^{{\tiny \alpha}}(\#7)$\quad\quad & \quad$S_{{\small L}%
}^{{\tiny \beta}}(\#10)$\quad\quad & \quad$S_{{\small L}}^{{\tiny \gamma}%
}(\#4)$\quad\quad\\\hline
\quad Type\quad\quad & [{\small 8 + 7}] & [{\small 2 + 13}] & [{\small 0 +
15}] & [{\small 8 + 7}]\\\hline
\end{tabular}

\end{center}

\begin{quotation}
Table 6. Four distinct (non-isomorphic) sedenion-type loops and their maximal
subloop compositions.\bigskip
\end{quotation}

\noindent From this we find that the loops $S_{{\small L}},$ $S_{{\small L}%
}^{{\tiny \alpha}}$, and $S_{{\small L}}^{{\tiny \beta}}$ differ in their
maximal subloop compositions. Hence they are not isomorphic. Note that the
loop $S_{{\small L}}^{{\tiny \gamma}}$ has the same subloop composition as
$S_{{\small L}}$ but analysis shows that they are not isomorphic. Therefore
all of these four subloops are distinct.

The sedenion-type loops generate 16-dimensional non-associative real
subalgebras of $\mathbb{T}$ with zero divisors [2]. This is due to the fact
that (up to isomorphism) each of them contains $\widetilde{O}_{{\small L}},$
and hence $\widetilde{\mathbb{O}}$ as a subalgebra. Such algebras (and their
subalgebras) have potential applications in both pure and applied mathematics
and in theoretical physics [7, 8].

\subsection{Properties of $\mathbf{T}_{{\protect\small L}}$ and its Subloops}

All Cayley-Dickson algebras of dimension $D\geq8$ are non-associative and
non-commutative. Analysis of the loop $T_{{\small L}}$ shows that it is a
NAFIL that satisfies the following identities (universally quantified
equations) listed in Table 7. Naturally, these identities are also satisfied
by all subloops of $T_{{\small L}}$ (both associative and
non-associative).\bigskip

\begin{center}
$%
\begin{tabular}
[c]{|c|c|c|}\hline
\textbf{Special Loop Property} & \textbf{Acronym} & \textbf{Defining
Equation}\\\hline\hline
\multicolumn{1}{|l|}{{\small Inverse Property}} & {\small IP} & $\ell
^{-1}(\ell q)=(q\ell)\ell^{-1}=q$\\\hline
\multicolumn{1}{|l|}{{\small Alternative Property}} & {\small AP} & $\ell(\ell
q)=\ell^{2}q${\small \ and }$(\ell q)q=\ell q^{2}$\\\hline
\multicolumn{1}{|l|}{{\small Flexible Law}} & {\small FL} & $\ell_{i}(\ell
_{k}\ell_{i})=(\ell_{i}\ell_{k})\ell_{i}$\\\hline
\multicolumn{1}{|l|}{{\small C Loop Property}} & {\small CP} & $\ell_{i}%
[\ell_{j}(\ell_{j}\ell_{k})]=[(\ell_{i}\ell_{j})\ell_{j}]\ell_{k}$\\\hline
\multicolumn{1}{|l|}{{\small Power Associative Property}} & {\small PAP} &
$\ell^{a}\cdot\ell^{b}=\ell^{a+b}$\\\hline
\multicolumn{1}{|l|}{{\small Weak Inverse Property}} & {\small WIP} &
$\ell(q\ell)^{-1}=q^{-1}$\\\hline
\multicolumn{1}{|l|}{{\small Anti-Automorphic Inverse Property}} &
{\small AAIP} & $(\ell q)^{-1}=q^{-1}\ell^{-1}$\\\hline
\end{tabular}
\ $\medskip
\end{center}

\begin{quotation}
Table 7. Some special loop properties (identities) satisfied by all
Cayley-Diction loops of order $m\geq16.$ The trigintaduonion loop $T_{L}$ and
its subloops satisfy these identities. The algebra $\mathbb{T},$ however,
satisfies only the flexible and power-associative identities. \bigskip
\end{quotation}

Not all of the identities satisfied by the loop $T_{{\small L}}$ listed in
Table 7 are satisfied by the algebra $\mathbb{T}$ as a whole. Moreover there
are identities satisfied only by its subloops. For instance, the octonion loop
and the quaternion and cyclic groups satisfy the Moufang identity.\bigskip

\subsection{Subloop Lattice of $T_{{\protect\small L}}$}

The subloop lattice of $T_{{\small L}}$ (and hence the basic \emph{subalgebra
lattice} of $\mathbb{T}$) is determined by its subloop composition. This
consists of its 373 subloops classified into 9 isomorphy classes:
$\mathcal{C}_{\cong}\{S_{{\small L}}\},$ $\mathcal{C}_{\cong}\{S_{{\small L}%
}^{\alpha}\},$ $\mathcal{C}_{\cong}\{S_{{\small L}}^{\beta}\},$ $\mathcal{C}%
_{\cong}\{S_{{\small L}}^{\gamma})\},$ $\mathcal{C}_{\cong}\{O_{{\small L}%
}\},$ $\mathcal{C}_{\cong}\{\widetilde{O}_{{\small L}}\},$ $\mathcal{C}%
_{\cong}\{Q_{{\small 8}}\},$ $\mathcal{C}_{\cong}\{C_{{\small 4}}\},$ and
$\mathcal{C}_{\cong}\{C_{{\small 2}}\}.$\bigskip

The diagram shown in Figure 1 indicates the general form of the subloop
lattice of $T_{{\small L}}$ in terms of the isomorphy classes of its subloops.
The determination of the complete lattice is a complicated problem that we are
now trying to address. So far, we know that all subloops of $T_{{\small L}}$
are normal. Such a loop has been shown to have a \emph{modular} lattice [6].

\begin{center}%
{\includegraphics[
natheight=3.908100in,
natwidth=2.788200in,
height=3.9081in,
width=2.7882in
]%
{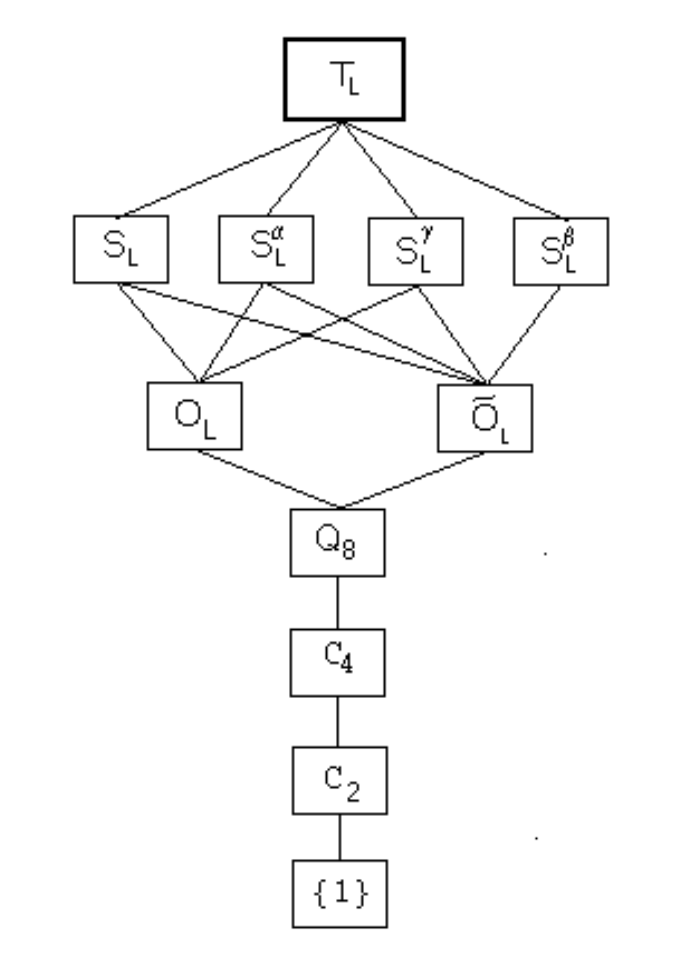}%
}%

\end{center}

\begin{quotation}
Figure 1. Lattice diagram (in block form) of the 9 isomorphy classes of the
subloops of the loop $T_{L}.\medskip$
\end{quotation}

\section{Summary}

The Cayley-Dickson algebra $\mathbb{T}$ of dimension $D=32$ (called the
\emph{trigintaduonions}) contains an embedded NAFIL loop $T_{{\small L}}$ of
order 64. This loop contains 31 subloops of order 32, 155 of order 16, 155 of
order 8, 31 of order 4, and 1 of order 2, all of which are normal. Such a loop
is known to have a modular lattice.

The subloops of orders 32 (\emph{sedenion-type}) and 16 (\emph{octonion-type})
are NAFILs, while those of orders 8, 4, and 2 are groups. These non-trivial
subloops of $T_{{\small L}}$ are classified into nine isomorphy
classes:\emph{\ }$\mathcal{C}_{\cong}\{S_{{\small L}}\},$\emph{\ }%
$\mathcal{C}_{\cong}\{S_{{\small L}}^{\alpha}\}$, $\mathcal{C}_{\cong%
}\{S_{{\small L}}^{\beta}\}$, $\mathcal{C}_{\cong}\{S_{{\small L}}^{\gamma}\}$
$\mathcal{C}_{\cong}\{O_{{\small L}}\},$ $\mathcal{C}_{\cong}\{\widetilde
{O}_{{\small L}}\},$ $\mathcal{C}_{\cong}\{Q_{{\small 8}}\},$ $\mathcal{C}%
_{\cong}\{C_{{\small 4}}\},$ and $\mathcal{C}_{\cong}\{C_{{\small 2}}\}.$

It can be shown that these subloops of $T_{{\small L}}$ of orders 32, 16, 8,
4, and 2 generate subalgebras of $\mathbb{T}$ of dimensions 16, 8, 4, 2, and
1, respectively. Hence the basic subalgebra structure of $\mathbb{T}$ is
directly related to the subloop structure of $T_{{\small L}}.$

Up to isomorphism, the octonion-type subloops of $T_{{\small L}}$ have the
same subloop composition: 7 quaternion groups $Q_{8}$, 7 cyclic groups $C_{4}%
$, and 1 cyclic group $C_{2}.$ Nevertheless, they differ in many fundamental
aspects (e.g. $O_{{\small L}}$.satisfies the Moufang identity while
$\widetilde{O}_{{\small L}}$ does not).

The 31 maximal subloops of order 32 form 4 isomorphy classes and 3 distinct
maximal subloop composition types: $[$8 $O_{{\small L}}$ + 7 $\widetilde
{O}_{{\small L}}]$, $[$2 $O_{{\small L}}$ + 13 $\widetilde{O}_{{\small L}}],$
and $[$0 $O_{{\small L}}$ + 15 $\widetilde{O}_{{\small L}}].$ This shows that
(up to isomorphism) every sedenion-type loop contains at least 7
quasi-octonion loops as subloops.

All 16-dimensional subalgebras of the trigintaduonions $\mathbb{T}$ (like
$\mathbb{S},$ $\mathbb{S}^{\alpha},$ $\mathbb{S}^{\beta},$and $\mathbb{S}%
^{\gamma}$) have zero divisors. Except for $\mathbb{S},$ these interesting
algebraic structures have not been previously identified. To date not much is
known about these algebras and it is therefore a challenge to determine their
structural features and their possible applications.

The lattice of the isomorphy classes of the subloops of $T_{L}$ is now known
(Figure 1). However, the complete subloop lattice of $T_{{\small L}}$ has not
yet been determined. This aspect of the problem is now being studied.\medskip

\begin{remark}
A recent review of the literature has shown that S. Catto and D. Chesley
(\emph{Twisted octonions and their symmetry groups}, Nuclear Physics
Proceedings Supplements, Volume 6, p. 428-432) have independently identified
the quasi-octonions which they have called the "twisted octonions." Their
approach to this problem is different from our method: it is based on the
analysis of quaternion triples. However, their results agree with our own
findings. Robert de Marrais (\emph{The 42 Assessors and the Box-Kites they
Fly: Diagonal axis-pairs systems of zero-divisors in the Sedenions},
http://arXiv.org/abs/math.GM/0011260) has also been studying the
Cayley-Dickson algebras using a different method which he calls the
\emph{Box-Kites} approach.
\end{remark}

{\footnotesize http://www.math.du.edu/loops/}\bigskip\bigskip\bigskip

{\footnotesize Alexander S. Carrascal, Lincoln A. Bautista, John P. Sta.
Maria, Jackie D. Urrutia, Bernadeth Nobles, College of Science, Polytechnic
University of the Philippines, Manila}\bigskip

\end{document}